# On a Kind of Nonlinear Diophantine Equation

*Massimo Salvi*

*Department of Mathematics and Computer Science*

*90123-Palermo, PA.*

*Italy*

e-mail: *massimo.salvi1@istruzione.it*

**Abstract.** In this paper we show a way to generalize the linear Diophantine equation $a_1 x_1 + a_2 x_2 + .. + a_n x_n = d$. We deal with the nonlinear Diophantine equation $\left|\begin{array}{c} A \\ X \end{array}\right| = \pm d$, which generalizes the linear one, and we give a necessary and sufficient condition for its solvabiliy. We show how the linear equation can be considered as a particular case of the nonlinear equation.

## 1. Introduction

It is known [1, p.30] that the diophantine equation:

$$a_1 x_1 + a_2 x_2 + .. + a_n x_n = d \qquad (1)$$

where $d$, $a_i$ and $x_i$ are required to be in Z, is solvable if and only if:

$$\gcd\{a_1, a_2 .., a_n\} \text{ divides } d \qquad (\text{Condition 1})$$

The generalization we want to introduce here starts from the observation that the equation $a_1 x_1 + a_2 x_2 = d$ can be written as $\left|\begin{array}{cc} a_1 & -a_2 \\ x_2 & x_1 \end{array}\right| = d$ (see [2] for other considerations), more generally we can consider equations of the form:

$$\left|\begin{array}{c} A \\ X \end{array}\right| = \pm d \qquad (2)$$

where $A$ and $X$ are integer matrices (see below). In this way we generalize the linear equation to a nonlinear equation. In this paper we show a necessary and sufficient condition for the solvability of this equation which generalizes the previous *Condition 1*. Moreover, we show that the possibility to write a linear diophantine equation (1) in the form (2) holds for equations with any possible number of unknowns (and not only two as in the example proposed). First we clear up some notations that we use in the paper.

Given any $n \times n$ matrix $M$, in this paper we denote by $|M|$ the determinant of $M$. In the case that $M$ is an integer and not square matrix, with $|M|$ we denote the greatest divisor (see below) of the matrix.

Be $A$ the matrix $r \times n$ $\begin{bmatrix} a_{11} & a_{12} & .. & a_{1n} \\ .. & & & \\ a_{r1} & a_{r2} & .. & a_{rn} \end{bmatrix}$ with n columns and $r$ rows, and $a_{ij} \in Z$, similarly be $X$ the matrix $\begin{bmatrix} x_{r+11} & x_{r+12} & .. & x_{r+1n} \\ .. & & & \\ x_{n1} & x_{n2} & .. & x_{nn} \end{bmatrix}$ with n columns and $n$-$r$ rows, and, $x_{ij} \in Z$

be $d > 0$ an integer, in this paper we deal with the Diophantine equation :

$$\left| \begin{matrix} A \\ X \end{matrix} \right| = \pm d$$

in which the entries $a_{ij}$ are supposed to be known integer values and the entries $x_{ij}$ are the unknowns to be found. We suppose to have, at least, a row of unknowns, hence $r < n$. In this paper we are interested in the entries $x_{ij}$ which satisfy $\left| \begin{matrix} A \\ X \end{matrix} \right| = d$ **or** $\left| \begin{matrix} A \\ X \end{matrix} \right| = -d$, for this reason we use the notation $\pm d$ in (2). We will refer to the matrix *X*, when it exists, as *solution* of the equation.

The equation (2) generalizes an equation studied by Smith [3] and before him by Hermite [4] in a more particular case. Smith analyzed the case in which *d*=1 and the greatest divisor of *A* (see below) is 1, Hermite analyzed the case in which *d*=1 and *A* is a vector of relatively prime integer numbers.

In this paper we give two main results:

1. A necessary and sufficient condition for the existence of integer solutions of (2);
2. We show that the equation (1) can be considered a particular case of the equation (2), and in this case the necessary and sufficient conditions for solvability of the equation (2) generalize the *Condition 1,* which holds for the equation (1).

In §2 we recall some concepts, definitions and theorems that are used in the rest of the paper, in §3 we give a necessary and sufficient condition for the solvability of equation (2), in §4 we show that the equation (2) actually generalizes the equation (1). We make use of various known results of matricial algebra (see [5], [6]).

## 2. Some Preliminaries Definitions and Results

In this section we recall some known results and definitions which are used in what follows.

### 2.1. Elementary Operations and Unimodular Matrices

Given an integer matrix *M*, an *elementary row operation* [7, p.12] is:

- Multiply a row by -1;
- Replace a row $R_i$ by $R_i + nR_j$, with $i \neq j$ and $n \in Z$ ; (3)
- Swap two rows.

The *elementary column operations* are defined in the same way by replacing the rows with the columns. Note, anyway, that the swapping of two rows (or two columns ), say $R_i$ and $R_j$, can be obtained by using only the first two operations by the sequence of operations: $R_i \to R_i + R_j$ ; $R_j \to R_i - R_j$ ; $R_i \to R_i - R_j$.

An integer matrix *U*, such that $|U| = \pm 1$, is called *unimodular matrix* ([8], [7, p.12]). It's easy to prove that the inverse of an unimodular matrix is unimodular, the product of two unimodular matrices is unimodular, the identity matrix I is clearly unimodular. That is, the set of unimodular $n \times n$ matrices, for a given *n*, provided with the matrices' multiplication operation is a group. Each of the above operations (3) is equivalent respectively to [7, p.13]:

- pre-multiplying *M* by an unimodular matrix *U* if we operate on the rows
- post-multiplying *M* by an unimodular matrix *U'* if we operate on the columns

A finite sequence of operations (3) on the rows or, respectively, on the colums of *M*, can be represented by an unimodular matrix, given by the product of all the matrices representing the elementary operations.

### 2.2. Lower Triangular Form of a Matrix

In §3 we prove Theorem 4, which concerns sufficient and necessary conditions for the solvability of equation (2). In order to do it we need the notion of *lower triangular form* of a matrix, hence we give the following [5] :

**Definition 1.** Given a full row rank integer matrix *M*, we say it is in Lower Triangular Form (LTF in what follows) if the matrix *M* is of the form $[N \quad O]$ where *O* is a matrix with all entries equal to 0 and *N* is a nonsingular square matrix of the form:

$$\begin{bmatrix} a_{11} & 0 & .. & 0 \\ & & 0 & \\ & & & 0 \\ a_{i1} & .. & & a_{ii} \end{bmatrix}$$

in which all the diagonal entries are positive and, for every row i, if $j > i$ we have $a_{ij} = 0$

**Example 1.** The following matrix is in LTF: $\begin{bmatrix} 1 & 0 & 0 & 0 & 0 \\ 0 & 3 & 0 & 0 & 0 \\ 1 & -1 & 2 & 0 & 0 \end{bmatrix}$

The following theorem can be proved by applying the algorithm which we show in §2.3 (see [9, cap.4] for other proofs):

**Theorem 1.** If $M$ is a full row rank integer matrix with $r$ rows and $c$ columns and $r \leq c$, then there exists an unimodular matrix $U$ such that the matrix $MU$ is in LTF.

## 2.3. Greatest Divisor of a Matrix

We give the following definitions [3]:

**Definition 2.** If $M$ is a matrix, the *determinants of M* are defined as the determinants of the greatest square matrices contained in $M$.

**Definition 3.** If $M$ is an integer matrix, the *greatest divisor* of $M$ is defined as the greatest common divisor of the determinants of $M$.

Note that if $M$ is a full row rank matrix then its greatest divisor is well defined and different from 0, if $M$ is not full row rank than all the determinants are 0, so the greatest divisor is not defined.

**Theorem 2.** If $M$ is a full row rank integer matrix and we pre-multiply or post-multiply $M$ by an unimodular matrix, the greatest divisor of the resulting matrix does not change.

See [3, p.311] for a proof.

If $M$ is in LTF, that is $M = [N \quad O]$, it's clear that the greatest divisor of $M$ can be obtained by the product of all the diagonal entries of the square matrix $N$. The reason for this is that $N$ is the only square matrix contained in $M$ with determinant not equal to 0. Hence, by using the previous Theorem 1 and Theorem 2 we have the:

**Theorem 3.** If $M$ is a full rank integer matrix with r rows and c columns and $r \leq c$, and $M' = [N \quad O]$ is the LTF of $M$, then the greatest divisor of $M$ is equal to the greatest divisor of $M'$, which is equal to the product of all the diagonal entries of $N$.

## 2.4. How to Reduce a Matrix into LTF

In this section we show a way to reduce a matrix into LTF by using the elementary operations (3) on the columns. We are not interested in computation efficiency, but only in showing a way to do it. The algorithm is based on the Euclidean Algorithm ([10], [9, p.52]). We explain and use the algorithm to reduce into lower triangular form a full row rank integer matrix M with r rows and c columns, and $r \leq c$.

Be $M=\begin{bmatrix} m_{11} & m_{12} & .. & .. & .. & m_{1c} \\ m_{21} & .. & .. & .. & .. & .. \\ .. & .. & .. & .. & .. & .. \\ m_{r1} & .. & .. & .. & .. & m_{rc} \end{bmatrix}$

If we consider the first row, all its entries cannot be all equal to 0 since the M is full rank, hence there exists the greatest common divisor $d_1$ of its c elements. We show that, by the following algorithm 1, $M$ can be reduced into the following form M' by unimodular operations on the columns.

$$M' = \begin{bmatrix} d_1 & 0 & .. & .. & .. & 0 \\ m'_{21} & .. & .. & .. & .. & .. \\ .. & .. & .. & .. & .. & .. \\ m'_{r1} & .. & .. & d_r & 0 & 0 \end{bmatrix} = \begin{bmatrix} N & 0 \end{bmatrix}$$

where $N$ is a lower triangular square matrix in which $d_i > 0, \forall i$.

Algorithm 1

First we repeat the following steps until in the first row there is only one nonzero entry, $d_1$:

1. Find the value of the row that has the smaller absolute value different from zero, and we denote it by $m_{1j}$;
2. For every $i \ne j$ by using the euclidean division algorithm we can write:
$m_{1i} = q_i m_{1j} + r_i$ where $0 \le r_i < |m_{1j}|$;
3. In each column $C_i$, $i \ne j$, we operate the substitution: $C_i \to C_i - q_i C_j$.

Then we can move $d_1$ in the first position of the row, and eventually change its sign, by elementary column operations. It can be proven that $d_1$ is the greatest common divisor of the entries of the first row [10]. Then we can continue this process starting from the matrix obtained at the end of the previous steps, without considering, this time, the first column and considering the entries of the second row. Then we can continue with the third row, and so on. Going on this way we can transform the matrix $M$ in a matrix of the form $M'$. All the operations are equivalent to post-multiply $M$ by an unimodular matrix $U$. Note that for all the diagonal entries we have $d_i > 0$, because if $d_i < 0$ we can multiply the column by -1, and for every $i$, $d_i \ne 0$ otherwise, for Theorem 2 and Theorem 3, the greatest divisor of $M$ cannot be defined and positive as required by the hypothesis that $M$ is full rank.

**Example 2.**

Let's consider the full row rank matrix $\begin{bmatrix} 2 & 2 & -3 & 4 \\ 2 & 2 & 1 & 2 \end{bmatrix}$ and apply the previous algorithm.

|  | Column transformation | Equivalent matrix post-multiplication |
|---|---|---|
|  | Algorithm 1 |  |
| T1 | $\begin{cases} C2 \to C2 - C1 \\ C3 \to -C3 - C1 \\ C4 \to C4 - 2C1 \end{cases}$ | $\begin{bmatrix} 2 & 2 & -3 & 4 \\ 2 & 2 & 1 & 2 \end{bmatrix} \cdot \begin{bmatrix} 1 & -1 & -1 & 2 \\ 0 & 1 & 0 & 0 \\ 0 & 0 & -1 & 0 \\ 0 & 0 & 0 & 1 \end{bmatrix} = \begin{bmatrix} 2 & 0 & 1 & 0 \\ 2 & 0 & -3 & -2 \end{bmatrix}$ |
| T2 | $\{C1 \to C1 - 2C3$ | $\begin{bmatrix} 2 & 0 & 1 & 0 \\ 2 & 0 & -3 & -2 \end{bmatrix} \cdot \begin{bmatrix} 1 & 0 & 0 & 0 \\ 0 & 1 & 0 & 0 \\ -2 & 0 & -1 & 0 \\ 0 & 0 & 0 & 1 \end{bmatrix} = \begin{bmatrix} 0 & 0 & 1 & 0 \\ 8 & 0 & -3 & -2 \end{bmatrix}$ |
| T3 | $\{C1 \leftrightarrow C3$ | $\begin{bmatrix} 0 & 0 & 1 & 0 \\ 8 & 0 & -3 & -2 \end{bmatrix} \cdot \begin{bmatrix} 0 & 0 & 1 & 0 \\ 0 & 1 & 0 & 0 \\ 1 & 0 & 0 & 0 \\ 0 & 0 & 0 & 1 \end{bmatrix} = \begin{bmatrix} 1 & 0 & 0 & 0 \\ -3 & 0 & 8 & -2 \end{bmatrix}$ |
| T4 | $\{C3 \to C3 + 4C4$ | $\begin{bmatrix} 1 & 0 & 0 & 0 \\ -3 & 0 & 8 & -2 \end{bmatrix} \cdot \begin{bmatrix} 1 & 0 & 0 & 0 \\ 0 & 1 & 0 & 0 \\ 0 & 0 & 1 & 0 \\ 0 & 0 & 4 & 1 \end{bmatrix} = \begin{bmatrix} 1 & 0 & 0 & 0 \\ -3 & 0 & 0 & -2 \end{bmatrix}$ |
| T5 | $\begin{cases} C2 \to -C4 \\ C4 \to C2 \end{cases}$ | $\begin{bmatrix} 1 & 0 & 0 & 0 \\ -3 & 0 & 0 & -2 \end{bmatrix} \cdot \begin{bmatrix} 1 & 0 & 0 & 0 \\ 0 & 0 & 0 & 1 \\ 0 & 0 & 1 & 0 \\ 0 & -1 & 0 & 0 \end{bmatrix} = \begin{bmatrix} 1 & 0 & 0 & 0 \\ -3 & 2 & 0 & 0 \end{bmatrix}$ |

Now, by multiplying all the unimodular matrices we obtain :

$$\begin{bmatrix} 2 & 2 & -3 & 4 \\ 2 & 2 & 1 & 2 \end{bmatrix} \cdot \begin{bmatrix} -1 & 2 & -5 & -1 \\ 0 & 0 & 0 & 1 \\ -1 & 0 & 2 & 0 \\ 0 & -1 & 4 & 0 \end{bmatrix} = \begin{bmatrix} 1 & 0 & 0 & 0 \\ -3 & 2 & 0 & 0 \end{bmatrix}$$

## 3. A Necessary and Sufficient Condition for Solvability of $\begin{vmatrix} A \\ X \end{vmatrix} = \pm d$

Let's consider the equation $\begin{vmatrix} A \\ X \end{vmatrix} = \pm d$, where $A$ is a full rank integer matrix with $r$ rows and $c$ columns with $r < c$ and $d$ is an integer greater than 0. We denote with $|A|$ the greatest divisor of $A$. (Note that, if $A$ is square, this notation coincides with the determinant of $A$, which, in this case, is the only determinant and therefore is the greatest divisor of the matrix). We prove the following:

**Theorem 4.** If the integer matrix $A$ is given, the equation (2) ($d > 0$) has integer solutions if and only if $A$ is full rank and $|A|$ divides $d$ (*Condition 2*).

*Proof.* We show first that the condition is sufficient, that is if $|A|$ divides $d$, then a solution exists.

Be $A'$ the LTF of $A$, then $A'$ is of the form $\begin{bmatrix} a_{11} & 0 & 0 & .. & 0 \\ .. & .. & 0 & .. & 0 \\ a_{r1} & & a_{rr} & .. & 0 \end{bmatrix}$

and, for Theorem 1, there exists $U$, unimodular, such that $A'=AU$. If $A$ is full rank then $|A|$ is defined and different from 0, therefore we can set $k = \dfrac{d}{|A|}$, which is supposed to be integer. Let's consider the matrix $B = \begin{bmatrix} 0 & N \end{bmatrix}$ with c columns and c-r rows, where N is a square diagonal matrix of the form:

$$\begin{bmatrix} 1 & 0 & .. & 0 \\ 0 & 1 & .. & .. \\ .. & .. & .. & 0 \\ 0 & .. & 0 & l \end{bmatrix}, \text{ in which the entries are } a_{ij} = \begin{cases} 1 & \text{if } i=j \neq c-r \\ l & \text{if } i=j=c-r \\ 0 & \text{otherwise} \end{cases} \text{ and } l = \begin{cases} k & \text{if } |U|=1 \\ -k & \text{if } |U|=-1 \end{cases} \quad (4)$$

Now consider the matrix $\begin{bmatrix} A' \\ B \end{bmatrix}$, this is a square $c \times c$ matrix, and its determinant is given by:

$$\left| \begin{matrix} A' \\ B \end{matrix} \right| = |A| \cdot l \quad (5)$$

thus if we consider the product with the unimodular matrix $U^{-1}$ we obtain:

$$\begin{bmatrix} A' \\ B \end{bmatrix} \cdot U^{-1} = \begin{bmatrix} A'U^{-1} \\ BU^{-1} \end{bmatrix} = \begin{bmatrix} A \\ BU^{-1} \end{bmatrix}$$

and, by computing the determinant of the product, we have:

$\left| \begin{bmatrix} A' \\ B \end{bmatrix} \cdot U^{-1} \right| = \left| \begin{matrix} A' \\ B \end{matrix} \right| \cdot |U^{-1}|$, for the previous (5) this is equal to $|A| \cdot l \cdot |U^{-1}|$ and for the definition

(4) it's equal to $|A| \cdot k$, but $|A| \cdot k = d$ so we can conclude that $\left| \begin{matrix} A \\ BU^{-1} \end{matrix} \right| = d$, therefore the matrix $\begin{bmatrix} BU^{-1} \end{bmatrix}$ is a solution of (2).

Now we prove that the condition of the theorem is necessary.

If the equation (2) has a solution, then $A$ has to be full rank, otherwise the determinant of the matrix $\begin{bmatrix} A \\ X \end{bmatrix}$ would be equal to 0. Let's suppose $A$ is full rank (hence its greatest divisor is defined) and $B$ is a solution of (2), then we have $\begin{vmatrix} A \\ B \end{vmatrix} = \pm d$. Since $d \neq 0$ for Theorem 1 there exists an unimodular matrix $U$ such that the product $\begin{bmatrix} A \\ B \end{bmatrix} \cdot U$ is in LTF. On one hand we have, considering that $U$ is unimodular and $|U| = \pm 1$, $\left| \begin{bmatrix} A \\ B \end{bmatrix} \cdot U \right| = \begin{vmatrix} A \\ B \end{vmatrix} \cdot |U| = \pm d$, on the other hand the determinant of a matrix in LTF is given by the product of all the entries on the diagonal. If we write more explicitly the product:

$$\begin{bmatrix} A \\ B \end{bmatrix} \cdot U = \begin{bmatrix} AU \\ BU \end{bmatrix} = \begin{bmatrix} A' \\ B' \end{bmatrix} = \begin{bmatrix} a'_{11} & 0 & .. & & .. & .. & 0 \\ .. & .. & 0 & & .. & .. & .. \\ .. & .. & a'_{rr} & & .. & .. & .. \\ .. & .. & b'_{r+1r} & b'_{r+1r+1} & .. & .. \\ .. & .. & .. & & .. & .. & 0 \\ b'_{n1} & .. & .. & & .. & .. & b'_{nn} \end{bmatrix}$$, then we can write:

$$\prod_{i=1}^{r} a'_{ii} \cdot \prod_{j=r+1}^{n} b'_{jj} = \pm d \qquad (6)$$

Now we have, from Theorem 3, that the product $\prod_{i=1}^{r} a'_{ii}$ is the greatest divisor of $A'$, which is equal, for Theorem 2, to the greatest divisor of $A$, denoted by $|A|$. Hence from (6) we can write $|A| \cdot \prod_{j=r+1}^{n} b'_{jj} = \pm d$, therefore $|A|$ divides $d$. ∎

**Example 3.**

Be the equation to solve: $\begin{vmatrix} 1 & 2 & -3 & 4 \\ 0 & 1 & 1 & 2 \\ x_{31} & x_{32} & x_{33} & x_{34} \\ x_{41} & x_{42} & x_{43} & x_{44} \end{vmatrix} = \pm 2$

The greatest divisor of the matrix $A$ is g.c.d $\{1,1,2,5,0,-10\} = 1$, and 1 divides 2. It means, by Theorem 4, that the equation is solvable.

**Observation 1.** It's easy to understand that the previous Theorem 4 holds for the equation $\begin{vmatrix} X \\ A \end{vmatrix} = \pm d$ as well. In fact the swapping of the matrix $A$ and $X$ within the determinant can be obtained by swapping the rows in a certain way, and this operations have the only effect to, eventually, change the sign.

**Observation 2.** Note that, if $A$ is not given the equation (which becames $|X|=\pm d$ ) is always solvable, even if $d=0$. In fact in this case we can always define, as we did in (4), a matrix which solve the equation $|X|=\pm d$:

$$\begin{bmatrix} 1 & 0 & .. & 0 \\ 0 & 1 & .. & .. \\ .. & .. & .. & 0 \\ 0 & .. & 0 & \pm d \end{bmatrix}, \text{ in which the entries are } a_{ij} = \begin{cases} 1 & \text{if } i=j \neq n \\ \pm d & \text{if } i=j=n \\ 0 & \text{otherwise} \end{cases}$$

## 4. The Equation (2) generalizes the Equation (1)

Let's go back to the equation (1) and let's address the problem of showing that both equation (2) and the conditions for its solvability *Condition 2*, found in the previous section, generalize the equation (1) and the *Condition 1*. To do this it's sufficient to prove that, given an equation of the form (1), there exists an integer matrix $A$ such that (1) can be written in the form (2). Namely, if we have:

$$a_1 x_1 + a_2 x_2 + .. + a_n x_n = \begin{vmatrix} a_{11} & .. & .. & .. & a_{1n} \\ . & & & & . \\ . & & & & . \\ a_{n-11} & .. & .. & .. & a_{n-1n} \\ x_1 & x_2 & .. & .. & x_n \end{vmatrix} = \begin{vmatrix} A \\ X \end{vmatrix}$$

then it's easy to understand that the *Condition 2* coincides with the *Condition 1*, in fact the set of determinants of $A$ coincides with the set $\{a_1, a_2.., a_n\}$, hence the gcd is equal. Moreover the previous equation entails that an equation of the form (1) can be considered as a particular equation of the form (2). We can prove the following:

**Theorem 5.** Given the linear form $a_1 x_1 + a_2 x_2 + .. + a_n x_n$, in which all values $\{a_1, a_2.., a_n\}$ are integers, it is always possible to find an integer matrix $A$, $(n-1) \times n$, such that:

$$a_1 x_1 + a_2 x_2 + .. + a_n x_n = \begin{vmatrix} a_{11} & .. & .. & .. & a_{1n} \\ . & & & & . \\ . & & & & . \\ a_{n-11} & .. & .. & .. & a_{n-1n} \\ x_1 & x_2 & .. & .. & x_n \end{vmatrix} = \begin{vmatrix} A \\ X \end{vmatrix} \qquad (7)$$

*Proof.* We can suppose that the $\gcd\{a_1, a_2.., a_n\}=1$, in fact given any form $a'_1 x_1 + a'_2 x_2 + .. + a'_n x_n$ and setting $k = \gcd\{a'_1, a'_2.., a'_n\}$ we can always consider the new form: $\frac{a'_1 x_1 + a'_2 x_2 + .. + a'_n x_n}{k} = a_1 x_1 + a_2 x_2 + .. + a_n x_n$ in which $\gcd\{a_1, a_2.., a_n\}=1$. If we suppose that a matrix $A$ exists such that the previous (7) holds, then by multiplying the first row (or any other) of the matrix $A$ by $k$, we obtain that (7) holds for $a'_1 x_1 + a'_2 x_2 + .. + a'_n x_n$ as well.

Now let's consider the matrix $M$ defined as follows:

$$M = \begin{bmatrix} 1 & 0 & .. & .. & 0 \\ 0 & 1 & 0 & .. & 0 \\ .. & 0 & .. & .. & 0 \\ 0 & .. & 0 & 1 & 0 \\ M_{n1} & M_{n2} & .. & .. & M_{nn} \end{bmatrix}$$

in the last row the entries are:
$$\begin{cases} M_{n1} = b_{11}x_1 + b_{12}x_2 + .. + b_{1n}x_n \\ M_{n2} = b_{21}x_1 + b_{22}x_2 + .. + b_{2n}x_n \\ .. \\ M_{nn-1} = b_{n-11}x_1 + b_{n-12}x_2 + .. + b_{n-1n}x_n \\ M_{nn} = a_1 x_1 + a_2 x_2 + .. + a_n x_n \end{cases}$$

where $b_{ij}$ are integer values which will be found in what follows,

the rest of the entries are: $\begin{cases} M_{ij} = 1 & \text{if } i = j \\ M_{ij} = 0 & \text{if } i \neq j \end{cases}$

By definition we have:
$$\begin{vmatrix} 1 & 0 & .. & .. & 0 \\ 0 & 1 & 0 & .. & 0 \\ .. & 0 & .. & .. & 0 \\ 0 & .. & 0 & 1 & 0 \\ M_{n1} & M_{n2} & .. & .. & M_{nn} \end{vmatrix} = M_{nn} = a_1 x_1 + a_2 x_2 + .. + a_n x_n$$

hence the theorem will be proved if we are able to find an unimodular matrix $U$:

$$U = \begin{bmatrix} u_{11} & u_{12} & .. & .. & u_{1n} \\ u_{21} & .. & .. & .. & .. \\ .. & .. & .. & .. & .. \\ .. & .. & .. & .. & .. \\ u_{n1} & u_{n2} & .. & .. & u_{nn} \end{bmatrix} \text{ such that } |U| = 1 \text{ and: } \begin{bmatrix} 1 & 0 & .. & .. & 0 \\ 0 & 1 & 0 & .. & 0 \\ .. & 0 & .. & .. & 0 \\ 0 & .. & 0 & 1 & 0 \\ M_{n1} & M_{n2} & .. & .. & M_{nn} \end{bmatrix} U = \begin{bmatrix} a_{11} & a_{12} & .. & .. & a_{1n} \\ a_{21} & .. & .. & .. & .. \\ .. & .. & .. & .. & .. \\ a_{n-11} & .. & .. & .. & a_{n-1n} \\ x_1 & x_2 & .. & .. & x_n \end{bmatrix}$$

as far as the last row is concerned, we have:

$$[b_1x_1+b_{12}x_2+..+b_{1n}x_n;\ \ ..\ \ ..\ \ b_{n-11}x_1+b_{n-12}x_2+..+b_{n-1n}x_n;\ a_1x_1+a_2x_2+..+a_nx_n]\begin{bmatrix} u_{11} & u_{12} & .. & .. & u_{1n} \\ u_{21} & .. & .. & .. & .. \\ .. & .. & .. & .. & .. \\ .. & .. & .. & .. & .. \\ u_{n1} & u_{n2} & .. & .. & u_{nn} \end{bmatrix}=[x_1\ \ x_2\ \ ..\ \ ..\ \ x_n]$$

by computing this product we obtain:

$$\begin{cases} (b_1u_{11}+b_2u_{21}+..+b_{n-1}u_{n-11}+a_1u_{n1})x_1+(b_{12}u_{11}+b_{22}u_{21}+..+b_{n-12}u_{n-11}+a_2u_{n1})x_2+..+(b_{1n}u_{11}+b_{2n}u_{21}+..+b_{n-1n}u_{n-11}+a_nu_{n1})x_n=x_1 \\ (b_1u_{12}+b_2u_{22}+..+b_{n-1}u_{n-12}+a_1u_{n2})x_1+(b_{12}u_{12}+b_{22}u_{22}+..+b_{n-12}u_{n-12}+a_2u_{n2})x_2+..+(b_{1n}u_{12}+b_{2n}u_{22}+..+b_{n-1n}u_{n-12}+a_nu_{n2})x_n=x_2 \\ .. \\ .. \\ (b_1u_{1n}+b_2u_{2n}+..+b_{n-1}u_{n-1n}+a_1u_{nn})x_1+(b_{12}u_{1n}+b_{22}u_{2n}+..+b_{n-12}u_{n-1n}+a_2u_{nn})x_2+..+(b_{1n}u_{1n}+b_{2n}u_{2n}+..+b_{n-1n}u_{n-1n}+a_nu_{nn})x_n=x_n \end{cases}$$

we can verify that the previous system is equivalent to requiring that:

$$\begin{bmatrix} u_{11} & u_{21} & .. & .. & u_{n1} \\ u_{12} & .. & .. & & u_{n2} \\ .. & & & & \\ .. & .. & & .. & \\ u_{1n} & & & & u_{nn} \end{bmatrix}\begin{bmatrix} b_{11} & b_{12} & & b_{1n} \\ b_{21} & b_{22} & & .. \\ & .. & .. & \\ b_{n-11} & b_{n-12} & & b_{n-1n} \\ a_1 & a_2 & & a_n \end{bmatrix}=\begin{bmatrix} 1 & 0 & .. & 0 & 0 \\ 0 & 1 & & & \\ .. & & .. & & 0 \\ .. & & & .. & 0 \\ 0 & & & 0 & 1 \end{bmatrix}=I \qquad (8)$$

or
$$U^T\begin{bmatrix} b_{11} & b_{12} & & b_{1n} \\ b_{21} & b_{22} & & .. \\ .. & .. & & \\ b_{n-11} & b_{n-12} & & b_{n-1n} \\ a_1 & a_2 & & a_n \end{bmatrix}=I$$

Now, since we required that $|U|=1$, it follows that $|U^T|=1$, and by using Binet's theorem in (8) we have:

$$\begin{vmatrix} b_{11} & b_{12} & & b_{1n} \\ b_{21} & b_{22} & & .. \\ .. & .. & & \\ b_{n-11} & b_{n-12} & & b_{n-1n} \\ a_1 & a_2 & & a_n \end{vmatrix}=1 \qquad (9)$$

But in the equation (9) $\gcd\{a_1,a_2..,a_n\}=1$, and 1 divides 1, so we can apply Theorem 4 (Observation 1) and conclude that (9) is solvable. Namely, there exists the required set of integer entries $b_{ij}$ and, therefore, there exists the required unimodular matrix $U$ as well. In fact from (8) we have:

$$(U^T)^{-1} = \begin{bmatrix} b_{11} & b_{12} & & b_{1n} \\ b_{21} & b_{22} & & .. \\ .. & .. & & \\ b_{n-11} & b_{n-12} & & b_{n-1n} \\ a_1 & a_2 & & a_n \end{bmatrix}$$

and, from what we have seen in 2.1, the inverse of an unimodular matrix is unimodular as well, hence $U^T$ is unimodular, therefore $U$ is unimodular as well. ∎

## 5. Concluding Remarks

In this paper we have shown that the linear diophantine equation: $a_1x_1 + a_2x_2 + .. + a_nx_n = d$ can be generalized to the nonlinear diophantine equation $\left| \begin{array}{c} A \\ X \end{array} \right| = \pm d$, where $A$ and $X$ are integer matrices and $d$ is an integer not equal to 0. We have shown that the necessary and sufficient condition for the solvability of the former, i.e. $\gcd\{a_1, a_2 .., a_n\}$ divides $d$, holds for the latter as well, if we use the notion of greatest divisor for an integer matrix.